\documentclass[12pt]{article}
\usepackage{amssymb, amsmath, amsthm, amscd}
\usepackage[dvips]{graphics}
\usepackage[utf8]{inputenc}
\usepackage[all,cmtip]{xy}
\usepackage{bbm}
\usepackage{enumitem}
\usepackage{setspace}
\usepackage[colorlinks=true]{hyperref}
\hypersetup{colorlinks   = true}
\hypersetup{linkcolor=blue}
\begin{document}

\newtheorem{The}{Theorem}[section]
\newtheorem{Lem}[The]{Lemma}
\newtheorem{Prop}[The]{Proposition}
\newtheorem{Cor}[The]{Corollary}
\newtheorem{Rem}[The]{Remark}
\newtheorem{Obs}[The]{Observation}
\newtheorem{SConj}[The]{Standard Conjecture}
\newtheorem{Titre}[The]{\!\!\!\! }
\newtheorem{Conj}[The]{Conjecture}
\newtheorem{Question}[The]{Question}
\newtheorem{Prob}[The]{Problem}
\newtheorem{Def}[The]{Definition}
\newtheorem{Not}[The]{Notation}
\newtheorem{Claim}[The]{Claim}
\newtheorem{Conc}[The]{Conclusion}
\newtheorem{Ex}[The]{Example}
\newtheorem{Fact}[The]{Fact}
\newtheorem{Formula}[The]{Formula}
\newtheorem{Formulae}[The]{Formulae}
\newcommand{\C}{\mathbb{C}}
\newcommand{\R}{\mathbb{R}}
\newcommand{\N}{\mathbb{N}}
\newcommand{\Z}{\mathbb{Z}}
\newcommand{\Q}{\mathbb{Q}}
\newcommand{\Proj}{\mathbb{P}}
\newcommand{\Rc}{\mathcal{R}}
\newcommand{\Oc}{\mathcal{O}}
\newcommand\houda[1]{{\textcolor{green}{#1}}}
\newcommand\dan[1]{{\textcolor{blue}{#1}}}
\begin{center}

{\Large\bf  Deformation Stability of $p$-SKT and $p$-HS manifolds}

\end{center}

\begin{center}

{\large Houda Bellitir}

\end{center}

\vspace{1ex}

\noindent{\small{\bf Abstract.} In this paper, we introduce the notions of $p$-Hermitian-symplectic and $p$-pluriclosed compact complex manifolds as generalisations for an arbitrary positive integer $p$ not exceeding the complex dimension of the manifold of the standard notions of Hermitian-symplectic and SKT manifolds that correspond to the case $p=1$. We then notice that these two properties are equivalent on $\partial\bar{\partial}$-manifolds and go on to prove that in (smooth) complex analytic families of $\partial\bar{\partial}$-manifolds, they are deformation open. Concerning closedness results, we prove that the cones $\mathcal{A}_p$, resp. $\mathcal{C}_p$, of Aeppli cohomology classes of strictly weakly positive $(p,p)$-forms $\Omega$ that are $p$-pluriclosed, resp. $p$-Hermitian-symplectic, must be equal on the limit fibre if they are equal on the other fibres and if some rather weak $\partial\bar{\partial}$-type assumptions are made on the other fibres.

\vspace{1ex}

\section{Introduction}

Let $\mathcal{X}$ be a complex manifold and let $\Delta$ be an open ball containing the origin in $\mathbb{C}^m$ for some $m\in\mathbb{N}^*$. Recall the following standard notion: a {\it complex analytic family of compact complex manifolds} is a proper holomorphic submersion $\pi: \mathcal{X}\longrightarrow \Delta$.

 This means that the following conditions are satisfied :
\begin{enumerate}
\item[(i)] For each $t\in\Delta$, $X_t:=\pi^{-1}(t)$ is a compact complex connected submanifold of $\mathcal{X}$.
\item[(ii)] the rank of Jacobian matrix of $\pi$ is equal to $n$ at every point of $\mathcal{X}$.
\item[(iii)] There is a locally finite open covering $\{\mathcal{U}_j / j=1,2,\cdots\}$ of $\mathcal{X}$ and complex-valued holomorphic functions $\xi_j^1(p),\cdots,\xi_j^n(p)$, defined on $\mathcal{U}_j$ such that for each $t$ , the set :
$$ \{p\in\mathcal{U}_j\longrightarrow(\xi_j^1(p),\cdots,\xi_j^n(p))/\mathcal{U}_j\cap\pi^{-1}(t)\neq\emptyset\}$$ 
is a system of local holomorphic coordinates of $X_t$.
\end{enumerate}
By a result of Ehresmann (\cite{Voi07}, Theorem 9.3), all the fibres $X_t:=\pi^{-1}(t)$, for all $t\in\Delta$, are $C^\infty$-diffeomorphic to a fixed $C^\infty$ manifold $X$. Therefore, the holomorphic family $(X_t)_{t\in\Delta}$ of compact complex manifolds can be viewed as a single $C^\infty$ manifold $X$ endowed with a $C^\infty$ family of complex structures $(J_t)_{t\in\Delta}$. \\

The purpose of this paper is to relax the notion of $p$-K\"ahlerianity introduced by Alessandrini and Andreatta in [AA87] for compact complex manifolds to two notions that offer analogues at the level of $(p,\,p)$-forms (satisfying a mild strict positivity assumption) for the classical notions of Hermitian-symplectic, respectively SKT (or pluriclosed) manifolds that correspond to the case $p=1$. Specifically, a smooth {\it strictly weakly positive} $(p,\,p)$-form $\Omega$ (cf. Definition \ref{SWP}) will be said to be {\it $p$-Hermitian-symplectic ($p$-HS)}, respectively {\it $p$-SKT} (or {\it $p$-pluriclosed)}, if $\Omega$ is the $(p,\,p)$-type component of a real $d$-closed $(2p)$-form, respectively if $\partial\bar\partial\Omega=0$. Compact complex manifolds carrying such forms will be termed {\it $p$-Hermitian-symplectic ($p$-HS) manifolds}, respectively {\it $p$-SKT (or $p$-pluriclosed) manifolds} in Definition \ref{def}.

\vspace{2ex}

Our first main result (cf. Theorems \ref{HS} and \ref{SKT} for precise statements) places these two notions in the context of deformations of complex structures.

\begin{The}\label{The:openness_introd} The notion of $p$-Hermitian-symplectic manifold and the notion of $p$-SKT $\partial\bar\partial$-manifold are both {\bf open} under holomorphic deformations of the complex structure.

\end{The}

Recall that the property of $p$-K\"ahlerianity of \cite{AA87} and \cite{AB91} is not deformation-open.

\vspace{3ex}

In the second part of the paper, we will investigate some deformation-closedness issues. To this end, we introduce the following positivity cones in the Aeppli cohomology of bidegree $(p,\,p)$ on any compact complex $n$-dimensional manifold $X$:

$$\mathcal{A}_p(X)=\{[\Omega]_A/ \Omega \mbox{ strictly weakly positive such that } \partial\bar{\partial}\Omega=0\}\subset H_A^{p,p}(X,\mathbb{R})\subset H_A^{p,p}(X,\mathbb{C}),$$ 
$$\mathcal{C}_p(X)=\{[\Omega]_A/ \Omega \mbox{ strictly weakly positive such that } \Omega \mbox{ is } p\mbox{-HS} \}\subset H_A^{p,p}(X,\mathbb{R})\subset H_A^{p,p}(X,\mathbb{C}).$$

\noindent The inclusion $\mathcal{C}_p(X)\subset\mathcal{A}_p(X)$ always holds trivially.

We prove (cf. Proposition \ref{Prop}) the following behaviour of these cones under holomorphic deformations of the complex structure. 

\begin{Prop}\label{Prop:cones_introd} Let $(X_t)_{t\in\Delta}$ be a holomorphic family of compact complex manifolds. If $\mathcal{A}_p(X_t) = \mathcal{C}_p(X_t)$ for all $t\in\Delta\setminus\{0\}$ and if some special cases of the $\partial\bar\partial$-Lemma (that we call the hypotheses $\tilde{H}_1,\cdots,\tilde{H}_{p+1}$) are satisfied for all $t\neq 0$ while the cohomology numbers $h_A^{p+k-1,p-k+1}(t):=\dim H_A^{p+k-1,p-k+1}(X_t,\mathbb{C}) $, $h_{BC}^{p+k,p-k+1}(t):=\dim H_{BC}^{p+k,p-k+1}(X_t,\mathbb{C})$ and $h_\partial^{p+k,p-k+1}(t):=\dim H_\partial^{p+k,p-k+1}(X_t,\mathbb{C})$ are independent of $t\in\Delta$ for all $k\in\{1,\cdots,p+1\}$, then $\mathcal{A}_p(X_0) = \mathcal{C}_p(X_0)$

\end{Prop}

\vspace{3ex}

We hope that these results will contribute to a possible resolution of Barlet's {\it Properness Conjecture} for the Barlet space of relative cycles associated with a holomorphic family of {\it class} $\mathcal{C}$ manifolds. Before presenting this conjecture let us recall that a compact complex manifold $X$ is said to be in the {\it (Fujiki) class $\mathcal{C}$} if there exists a proper holomorphic bimeromorphic map (i.e. a modification) $\mu : \tilde{X}\longrightarrow X$ from a compact K\"ahler manifold $\tilde{X}$ to $X$ (see \cite{Var86}).\\

Fujiki introduced class $\mathcal{C}$ manifolds $X$, in \cite{Fuj78}, as meromorphic images of compact K\"ahler manifolds. This definition was proved to be equivalent to the above-stated definition by Varouchas in \cite{Var86}.\\
On the other hand, denote by $\mathcal{C}^p(\mathcal{X}\diagup\Delta)$ the relative Barlet space of effective analytic $p$-cycles contained in the fibres $X_t$, and let: 
$$\mathcal{C}(\mathcal{X}\diagup\Delta)=\underset{{0\leq p \leq n}}{\cup} \mathcal{C}^p(\mathcal{X}\diagup\Delta) $$ 
\begin{Conj}
(\cite{Pop10}, one of the versions of Barlet's conjecture) Let $\pi : \mathcal{X}\longrightarrow \Delta$ be a complex analytic family of compact complex manifolds such that the fibre $X_t$ is a class $\mathcal{C}$ manifold for every $t\in\Delta$. Then the irreducible components of the relative Barlet space $\mathcal{C}(\mathcal{X}\diagup\Delta)$ of cycles on $\mathcal{X}$ are proper over $\Delta$ in the following sense. Consider the holomorphic map
$$ \begin{array}{ll}
P :  &\mathcal{C}(\mathcal{X}/\Delta)\longrightarrow \Delta \\
   &  \hspace*{1.3cm} Z_t \mapsto P(Z_t)=t
\end{array}$$
mapping every divisor $Z_t\subset X_t$ contained in some fibre $X_t$ to the base point $t\in \Delta$. The map $P$ has the property that its restrictions to the irreducible components of $\mathcal{C}(\mathcal{X}/\Delta)$ are proper.
\end{Conj}

 In the end, we prove (cf. Corollary \ref{Cor}) that on a fixed compact complex manifold $X$ satisfying some assumptions weaker than the $\partial\bar{\partial}$-lemma, the natural map:
 $$H_{DR}^k(X,\mathbb{C})\oplus H_{DR}^{2n-k}(X,\mathbb{C})\longrightarrow \bigoplus_{p+q=k}H_A^{p,q}(X,\mathbb{C})\oplus\bigoplus_{p+q=2n-k}H_A^{p,q}(X,\mathbb{C})$$
is an {\it isomorphism}.\\

\noindent {\bf Acknowledgements}. The author is grateful to her supervisor Dan Popovici for his continuous guidance, suggestions, support, encouragements during her thesis and for a very careful reading of the paper as well as to Professor Ahmed Zeriahi for many interesting discussions.

 \section{Preliminaries}

In this section, we recall some background that will be needed in the sequel and introduce our first definitions. Let $X$ be a compact complex manifold of complex dimension $n$. 

The following positivity notion is standard, though a different name has been used in the literature, and goes back to Lelong. We spell out the details of the definition in order to dispel the confusion over the English and French meanings of ``positivity''.
\begin{Def}\label{SWP} (see e.g. \cite{Dem}, Chapter III and \cite{AB91}, Definition 1.1)
 (i)\, Let $V$ be a complex vector space of dimension $n$ and let $V^*$ be its dual. Fix any integer $1\leq p\leq n-1$. A $(p,p)-$form $\alpha\in \Lambda^{p,p}V^*$ is said to be {\bf strictly weakly positive} if for all linearly independent $\tau_1,\dots , \tau_q\in V^*$, with $q=n-p$, the $(n,n)-$form $$\alpha\wedge i \tau_1\wedge \bar{\tau_1}\wedge \cdots \wedge i\tau_q\wedge\bar{\tau}_q $$
is positive (i.e. $>0$) whenever $ \alpha\wedge i \tau_j\wedge \overline{\tau_j}\neq 0$ for all $j\in \{1,\cdots,q\}$.

\vspace{1ex}

(ii)\, Let $X$ be a complex manifold of dimension $n$. A smooth $(p,\,p)$-form $\alpha\in C^\infty_{p,\,p}(X,\,\C)$ on $X$ is said to be {\bf strictly weakly positive} if, for every point $x\in X$, $\alpha(x)\in \Lambda^{p,p}T_x^*X$ is a strictly weakly positive $(p,p)$-form on the holomorphic tangent space $T_xX$ to $X$ at $x$. 

\vspace{1ex}

(These forms are called {\bf transverse} in \cite{AB91}.)

\end{Def}

 It is standard that a $(p,\,p)$-form $\alpha\in \Lambda^{p,p}V^*$ on a vector space $V$ is strictly weakly positive if and only if its restriction $\alpha_{|E}$ to every $p$-dimensional vector subspace $E\subset V$ is a positive (i.e. $>0$) volume form on $E$. (See e.g. \cite{Dem}, Chapter III).

Consequently, a smooth $(p,\,p)$-form $\alpha\in C^\infty_{p,\,p}(X,\,\C)$ on a manifold $X$ is strictly weakly positive if and only if, for every coordinate patch $U\subset X$ and every $p$-dimensional complex submanifold $Y\subset U$, its restriction $\alpha_{|Y}$ is a positive (i.e. $>0$) volume form on $Y$. (See e.g. \cite{AB91}).

It is also standard that, for $p=1$ and $p=n-1$, the notion of strict weak positivity coincides with the usual notion of positive definiteness for $(1,\,1)$ and $(n-1,\, n-1)$-forms ($=$ the positivity of all the eigenvalues of the coefficient matrix). Moreover, all strictly weakly positive $(p,\,p)$-forms $\alpha$ are always {\it real}, in the sense that $\alpha = \overline\alpha$. (See e.g. \cite{Dem}, Chapter III).

\vspace{3ex}

 Now, we introduce $p$-Hermitian-symplectic ($p$-HS) and $p$-SKT forms and manifolds as the following generalisations of the notions of $p$-K\"ahler forms and manifolds introduced in \cite{AA87} and further studied in \cite{AB91}.

\begin{Def}\label{def} Let $X$ be a compact complex manifold of complex dimension $n$ and let $\Omega$ be a $C^\infty$ {\bf strictly weakly positive $(p,p)-$form} on $X$.
\begin{enumerate}
\item[(i)] $\Omega$ is said to be a {\bf $p-$Hermitian-symplectic form} if there exist forms $\alpha^{i,2p-i}\in C_{i,2p-i}^\infty(X,\mathbb{C})$ for $i\in\{0,\cdots,p-1\}$ such that $d(\sum_{i=0}^{p-1}\alpha^{i,2p-i}+\Omega+\sum_{i=0}^{p-1}\overline{\alpha^{i,2p-i}})=0$.\\
\item[(ii)] $\Omega$ is said to be a {\bf $p-$SKT form} if $\partial\bar{\partial}\Omega=0$.\\
\item[(iii)]$X$ is said to be a {\bf $p-$SKT manifold} (resp. a {\bf $p-$HS manifold}) if there exists a $p$-SKT (resp. a $p$-HS) strictly weakly positive $(p,p)$-form on $X$.\\
\item[(iv)] (\cite{AB91}) A compact complex manifold $X$ is said to be {\bf $p$-K\"ahler} if it supports a $d$-closed strictly weakly positive $(p,p)-$form $\Omega$. Such an  $\Omega$ is called $p$-{\it K\"ahler form}.
\end{enumerate}
\end{Def}

\vspace{3ex}

Finally, recall the following standard definitions. For all $p,q=0,\cdots,n$, the {\bf Bott-Chern cohomology group} of $X$ of type $(p,q)$ is defined as:
$$ H^{p,q}_{BC}(X,\mathbb{C})=\dfrac{\ker\{\partial: C^\infty_{p,q}(X)\longrightarrow C^\infty_{p+1,q}(X) \}\cap \ker\{\bar{\partial}: C^\infty_{p,q}(X)\longrightarrow C^\infty_{p,q+1}(X)\}}{Im\{\partial\bar{\partial}: C^\infty_{p-1,q-1}(X,\mathbb{C})\longrightarrow C^\infty_{p,q}(X,\mathbb{C})\}},$$
while the {\bf Aeppli cohomology group} of type $(p,q)$ is defined as:
$$  H^{p,q}_{A}(X,\mathbb{C})=\dfrac{\ker\{\partial\bar{\partial}: C^\infty_{p,q}(X,\mathbb{C})\longrightarrow C^\infty_{p+1,q+1}(X,\mathbb{C})\}}{Im\{\partial: C^\infty_{p-1,q}(X)\longrightarrow C^\infty_{p,q}(X) \}+ Im\{\bar{\partial}: C^\infty_{p,q-1}(X)\longrightarrow C^\infty_{p,q}(X)\}}.$$

\section{Deformations of pluriclosed manifolds}

Let us now recall some basic notions in the following

\begin{Def} Let $\omega>0$ be a Hermitian metric on a complex manifold $X$.
\begin{enumerate}
\item[1-]$\omega$ is said to be {\bf pluriclosed} (or {\bf SKT}) if $\partial\bar{\partial}\omega=0$.
\item[2-]$\omega$ is called {\bf Hermitian-symplectic} ({\bf H-S} for simplicity) (cf. Definition in \cite{ST10}) if there exists $\alpha^{0,2}\in C_{0,2}^\infty(X,\mathbb{C})$ such that  $$d(\overline{\alpha^{0,2}}+\omega+\alpha^{0,2})=0.$$ This is equivalent to the existence of a $C^\infty$ real $2$-form $ \Omega$ on $X$ such that $d\Omega=0$ and $\Omega^{1,1}=\omega>0$.
\item[3-] $X$ is a {\bf SKT manifold} (resp. {\bf H-S manifold}) if there exists an SKT (resp. a H-S) metric on $X$.
\item[4-] $X$ is a {\bf $\partial\bar{\partial}$-manifold} (cf. e.g. \cite{Pop14} for this version of a notion introduced in \cite{DGMS75}) if the $\partial\bar{\partial}$-lemma is satisfied on $X$, that is for all $ p$, $q$ and for every  $C^\infty$ $(p,q)$-form $u$ on $X$ such that $du=0$, the following properties are equivalent:
 $$\mbox{u is } d-exact \Leftrightarrow \mbox{u is } \partial-exact \Leftrightarrow \mbox{u is }  \bar{\partial}-exact \Leftrightarrow \mbox{u is } \partial\bar{\partial} -exact $$ 
\end{enumerate}

\end{Def}

It is proved in \cite{Pop15}, that on every $\partial\bar{\partial}$-manifold, the notion of SKT metric and the notion of Hermitian-symplectic metric are equivalent: 

\begin{Lem}\label{Lem:SKT_HS_equiv_ddbar}
Let $X$ be a  $\partial\bar{\partial}-$manifold. For any Hermitian metric $\omega$ on $X$, the following equivalence holds:\\
$$\omega \mbox{ is SKT }  \Longleftrightarrow \, \omega \mbox{ is Hermitian-symplectic}.$$
\end{Lem}

\noindent \noindent {\it Proof}. $\Leftarrow)$ This implication holds on any compact complex manifold.\\
Suppose that there exists $\alpha^{0,2}\in C_{0,2}^\infty(X,\mathbb{C})$  such that
 $d(\alpha^{0,2}+\omega+\overline{\alpha^{0,2}})=0$. This is equivalent to
$$
\left\{
\begin{array}{ll}
 \partial\overline{\alpha^{0,2}}=0 \\
 \partial\alpha^{0,2}+\bar{\partial}\omega=0 .
\end{array}
\right.
$$
Applying $\partial$ on the above identity we get $\partial\bar{\partial}\omega=0$.\\

\noindent Conversely, suppose that  $\partial\bar{\partial}\omega=0$ and $X$ is a $\partial\bar{\partial}-$manifold, then $\partial \omega \in \ker \bar\partial$.  Meanwhile $\partial \omega$ is a  $(2,1)$-form which is $d$-closed and $\partial$-exact, hence by the $\partial\bar{\partial}$-lemma $\partial \omega$ must be $\bar{\partial}-$exact.\\
This means that there exists $\alpha^{2,0}\in C_{2,0}^\infty(X,\mathbb{C})$  such that $\partial\omega=-\bar{\partial}\alpha^{2,0}=-\overline{\partial\alpha^{0,2}}$ with $\alpha^{0,2}:=\overline{\alpha^{2,0}}$. We get:
\begin{align*}
d(\alpha^{0,2}+\omega+\overline{\alpha^{0,2}})&=\partial\alpha^{0,2}+\bar{\partial}\alpha^{0,2}+\partial\omega+\bar{\partial}\omega+\partial\overline{\alpha^{0,2}}+\bar{\partial}\overline{\alpha^{0,2}}\\
&=(\partial\omega+\bar{\partial}\overline{\alpha^{0,2}})+(\overline{\partial\omega+\bar{\partial}\overline{\alpha^{0,2}}})+\bar{\partial}\alpha^{0,2}+\partial\overline{\alpha^{0,2}}\\
&=\bar{\partial}\alpha^{0,2}+\partial\overline{\alpha^{0,2}}.
\end{align*}
It remains to show that $\partial\alpha^{2,0}=0$. Indeed, we have $\bar{\partial}\partial\alpha^{2,0}=-\partial\bar{\partial}\alpha^{2,0}=\partial^2\omega=0$. So $\partial\alpha^{2,0}$ is a $d$-closed and $\partial$-exact $(3,0)$-form, therefore by the $\partial\bar{\partial}$-lemma it is $\bar{\partial}$-exact. For bidegree reasons, this means that  $\partial\alpha^{2,0}=0$. This proves that $\omega$ is H-S. \hfill $\Box$

\vspace{2ex}

In the rest of this section, we will show that the property of a manifold being Hermitian-symplectic is open under holomorphic deformations. 

\begin{The}\label{The:HS_open} Let $(X_t)_{t\in\Delta}$ be a holomorphic family of compact complex manifolds over an open disc $\Delta\subset\C$ containing the origin.

If there exists an H-S metric $\omega_0$ on $X_0$, then after possibly shrinking $\Delta$ about $0$, there exists a $C^\infty$ family $(\omega_t)_{t\in\Delta}$ of H-S metrics on the fibres $(X_t)_{t\in\Delta}$ whose term corresponding to $t=0$ is the original $\omega_0$. 

\end{The}

\noindent \noindent {\it Proof}.
Suppose that there exists a H-S metric $\omega_0$ on $X_0$. Then, there exists $ \alpha_0^{0,2}\in C_{0,2}^\infty(X_0,\mathbb{C}) $ such that $$d(\alpha_0^{0,2}+\omega_0+\overline{\alpha_0^{0,2}})=0$$
Putting $\Omega:=\alpha_0^{0,2}+\omega_0+\overline{\alpha_0^{0,2}}$, we get $\Omega$ a $d-$closed real $2-$form on $X_0$. Since, by Ehresmann, $\Omega$ does not depend on the complex structure of the fibre. Let $(\Omega_t^{1,1})_{t\in\Delta}$ be the $C^\infty$ family of components of $\Omega$ of $J_t$-type $(1,1)$, namely  $\Omega_t^{1,1}=(\Omega)_t^{1,1}$, so $(\Omega_t^{1,1})_{t\in\Delta}$ vary smoothly with $t$. Thus, by the continuity of the family $(\Omega_t^{1,1})_{t\in\Delta}$, the strict positivity of $\Omega_0^{1,1}$ implies the strict positivity of $\Omega_t^{1,1}$ for all $t\in\Delta$ sufficiently close to 0. So, we obtain a H-S metric $\omega_t$ on $X_t$ for every $t\in\Delta$ close to $0$. \hfill $\Box$

\vspace{2ex}

 Now, recall the following theorem of Wu \cite{Wu} according to which the $\partial\bar{\partial}$-property of compact complex manifolds is open under holomorphic deformations of the complex structure:

$$X_0 \mbox{ is } \partial\bar{\partial}\mbox{-manifold }  \Longrightarrow \, X_t \mbox{ is }\partial\bar{\partial}\mbox{-manifold for } t\in\Delta,\, t\thicksim 0.$$

\vspace{3ex}

\noindent As a consequence of Wu's theorem and of our Theorem \ref{The:HS_open}, we get the following 

\begin{Cor} Let $(X_t)_{t\in\Delta}$ be a holomorphic family of compact complex manifolds over an open disc $\Delta\subset\C$ containing the origin.

If $X_0$ is an SKT $\partial\bar{\partial}$-manifold, then $X_t$ is an SKT $\partial\bar{\partial}$-manifold for all $t\in\Delta$ sufficiently close to $0$. 

Moreover, if $X_0$ is an SKT $\partial\bar{\partial}$-manifold, any SKT metric $\omega_0$ on $X_0$ deforms in a $C^\infty$ way to a family $(\omega_t)_t$ of SKT metrics on the nearby fibres $(X_t)_t$.
\end{Cor}

\noindent {\it Proof.} On $\partial\bar{\partial}$-manifolds, the notions of H-S and SKT metrics are equivalent, as noticed in Lemma \ref{Lem:SKT_HS_equiv_ddbar}.  \hfill $\Box$

\section{Deformations of $p$-HS and $p$-SKT manifolds }
In this section, we will prove that the property of a compact complex $\partial\bar{\partial}$-manifold carrying $p$-SKT form is deformation open.

\begin{Lem}
Let $X$ be a $\partial\bar{\partial}-$manifold with $\dim_{\mathbb{C}}X=n$. Fix $p\in\{1,\cdots,n-1\}$.\\
If there exists a $p-$SKT form $\omega$ on $X$, then $\omega$ is a $p$-Hermitian-symplectic form.
\end{Lem}

\noindent \noindent {\it Proof}.
Suppose there exists a $p$-pluriclosed form $\omega$ on $X$. Namely, there exists a $C^\infty$ strictly weakly positive real $(p,p)-$form  $\omega$ such that $\partial\bar{\partial}\omega=0$. Notice that $\partial\omega$ is a $d$-closed $(p+1,p)-$form. By the $\partial\bar{\partial}-$lemma, $\partial\omega$ is $\bar{\partial}-$exact, that is there exists $\alpha^{p+1,p-1}\in C_{p+1,p-1}^\infty(X,\mathbb{C})$ such that  $\partial\omega=-\bar{\partial}\alpha^{p+1,p-1}$. Note that $\partial\alpha^{p+1,p-1}$ is $d-$closed, so by $\partial\bar{\partial}-$lemma there exists
$\alpha^{p+2,p-2}\in C_{p+2,p-2}^\infty(X,\mathbb{C})$ such that $\partial\alpha^{p+1,p-1}=-\bar{\partial}\alpha^{p+2,p-2}$.\\

We continue until we find a form $\alpha^{2p,0}\in C_{2p,0}^\infty(X,\mathbb{C}) $  such that  $\partial\alpha^{2p-1,1}=-\bar{\partial}\alpha^{2p,0}$. Thus, $\partial\alpha^{2p,0}$ is $d$-closed  hence by $\partial\bar{\partial}-$lemma it must also be $\bar{\partial}-$exact. However, the only $\bar{\partial}-$exact $(2p,0)-$form is zero, thus $\partial\alpha^{2p,0}=0$. As a consequence, we have:\\
$$
\begin{array}{cccc}
&d(\sum_{i=0}^{p-1}\alpha^{i,2p-i}+\omega+\sum_{i=0}^{p-1}\overline{\alpha^{i,2p-i}})\\
&=\partial\alpha^{0,2p}+\bar{\partial}\alpha^{0,2p}+\partial\overline{\alpha^{0,2p}}+\bar{\partial}\overline{\alpha^{0,2p}}+\partial\alpha^{1,2p-1}+\bar{\partial}\alpha^{1,2p-1}+\partial\overline{\alpha^{1,2p-1}}\\
&+\bar{\partial} \overline{\alpha^{1,2p-1}}+\cdots+\partial\alpha^{p-3,p+3}+\bar{\partial}\alpha^{p-3,p+3}+\partial\overline{\alpha^{p-3,p+3}}+\bar{\partial} \overline{\alpha^{p-3,p+3}}\\
&+\partial\alpha^{p-2,p+2}+\bar{\partial}\alpha^{p-2,p+2}+\partial\overline{\alpha^{p-2,p+2}}+\bar{\partial} \overline{\alpha^{p-2,p+2}}+\bar{\partial}\alpha^{p-1,p+1}\\
&+\partial\alpha^{p+1,p-1}=0.
\end{array}
$$

\noindent Thus, $\omega$ is a $p$-HS form.  \hfill $\Box$

\vspace{2ex}

A key observation for us will be the following deformation-openness property of the strict weak positivity for $(p,\,p)$-forms.

\begin{Lem}\label{positivity_def-open} Let $(X_t)_{t\in\Delta}$ be a holomorphic family of compact complex $n$-dimensional manifolds over an open disc $\Delta\subset\C$ containing the origin. Fix any $p\in\{1,\dots , n-1\}$. Let $(\Omega_t)_{t\in\Delta}$ be a $C^\infty$ family of real $(2p)$-forms on the fibres $X_t$ such that every $\Omega_t$ is of type $(p,\,p)$ for the complex structure $J_t$ of $X_t$.  

If $\Omega_0$ is strictly weakly positive for the complex structure $J_0$ of $X_0$, then $\Omega_t$ is strictly weakly positive for the complex structure $J_t$ of $X_t$ for all $t\in\Delta$ sufficiently close to $0\in\Delta$.

\end{Lem}

\noindent {\it Proof.} Let ${\cal X}=\cup_{t\in\Delta}X_t$ be the total space of the family and let ${\cal X}= U^{(1)}\cup\dots\cup U^{(N)}$ be an open covering by coordinate patches of ${\cal X}$ such that $X_0\subset U^{(1)}\cup\dots\cup U^{(N)}$. For every $\nu=1,\dots , N$, let $(z_1^{(\nu)},\cdots, z_n^{(\nu)}, \,t)$ be a local holomorphic coordinate system on $U^{(\nu)}\subset{\cal X}$, where $t$ is a holomorphic coordinate on a neighbourhood of $0$ in $\Delta$. Put $U^{(\nu)}_t:=U^{(\nu)}\cap X_t\subset X_t$ and $z_j^{(\nu)}(t):=z_{j|X_t}^{(\nu)}$.

Now, the strict weak positivity for $(p,\,p)$-forms on a complex manifold is a pointwise, hence local, property. For a fixed $t'\in\Delta$ and a fixed index $\nu\in\{1,\dots , N\}$, let $\tau_1(t'),\dots , \tau_q(t')$ (where $q=n-p$) be lineraly independent $C^\infty$ $J_{t'}$-$(1,\,0)$-forms on $U^{(\nu)}_{t'}\subset X_{t'}$. We extend the forms $\tau_j(t')$ to $C^\infty$ $(1,\,0)$-forms $\tau_1,\dots , \tau_q$ on $U^{(\nu)}\subset{\cal X}$.  (This is always possible for forms defined on a coordinate patch. For example, we can extend their coefficients by requiring them to remain constant in the $t$-direction.) For every $t\in\Delta$, the restrictions $\tau_1(t),\dots , \tau_q(t)$ of the forms $\tau_1,\dots , \tau_q$ to $U^{(\nu)}_t\subset X_t$ are of type $(1,\,0)$ for $J_t$. They are also linearly independent if $t$ lies in a sufficiently small neighbourhood (depending on the forms $\tau_1(t'),\dots , \tau_q(t')$) of $0$. Since $\Omega_0$ is strictly weakly positive for $J_0$, we have

$$\Omega_0\wedge i\tau_1(0)\wedge \overline{\tau_1}(0)\wedge\cdots\wedge i \tau_{n-p}(0)\wedge\overline{\tau_{n-p}}(0)>0 \hspace{3ex} \mbox{on} \hspace{1ex} U^{(\nu)}_0\subset X_0.$$

\noindent By continuity, we still have $\Omega_t\wedge i\tau_1(t)\wedge \overline{\tau_1}(t)\wedge\cdots\wedge i \tau_{n-p}(t)\wedge\overline{\tau_{n-p}}(t)>0$ for all $t$ lying in a sufficiently small neighbourhood (depending on the forms $\tau_1(t'),\dots , \tau_q(t')$) of $0$.

Meanwhile, the strict weak positivity at a point has to be tested on all the $p$-dimensional complex vector subspaces of an $n$-dimensional space. These subspaces form a {\it compact} complex manifold, the Grassmannian $G_{p,\,n}$. Thanks to the compactness of this Grassmannian, there exists a {\it uniform} open neighbourhood $\Delta_0$ of $0$ in $\Delta$ such that $\Omega_t$ is strictly weakly positive for the complex structure $J_t$ of $X_t$ for all $t\in\Delta_0$.  \hfill $\Box$

\vspace{3ex}

\noindent The first main result of this section is the following
\begin{The}\label{HS}
Let $X$ be a compact complex manifold of dimension $n$. Fix $p\in \{1,\cdots,n-1\}$. If there exists a $p$-HS form $\omega_0$ on $X_0$, then there exists a $p$-HS form $\omega_t$ on $X_t$ for $t\in\Delta$ sufficiently close to $0$.
\end{The}

\noindent \noindent {\it Proof}.
Assume that there exists a $p-$HS form $\omega_0$ on $X_0$. This means that there exist $\alpha^{i,2p-i}\in C_{i,2p-i}^\infty(X_0,\mathbb{C})$ for $i=0,\cdots,2p$ such that  $d(\sum_{i=0}^{p-1}\alpha^{i,2p-i}+\omega_0+\sum_{i=0}^{p-1}\overline{\alpha^{i,2p-i}})=0$. 

Put $\Omega:=\sum_{i=0}^{p-1}\alpha^{i,2p-i}+\omega+\sum_{i=0}^{p-1}\overline{\alpha^{i,2p-i}}$. This is a $d$-closed real $(2p)$-form on the $C^\infty$ manifold $X$ underlying the fibres $X_t$. Let $(\Omega_t^{p,p})_{t\in\Delta}$ be the $C^\infty$ family of components of $\Omega$ of $J_t$-type $(p,p)$. The forms $\Omega_t^{p,p}$ vary smoothly with $t\in\Delta$. 

Thanks to Lemma \ref{positivity_def-open}, the weak strict positivity of $\Omega^{p,p}_0=\omega_0$ implies the weak strict positivity of $\Omega_t^{p,p}$ for all $t\in\Delta$ sufficiently close to 0. Thus, $\omega_t:=\Omega_t^{p,p}$ is a $p$-HS form on $X_t$ for all $t\in\Delta$ sufficiently close to 0.   \hfill $\Box$

\vspace{2ex}

\noindent As a consequence we get
\begin{Conc}
If $X_0$ is a $p$-SKT $\partial\bar{\partial}$-manifold, then $X_0$ is a $p$-HS manifold. Therefore, $X_t$ is a $p$-HS manifold, hence also a $p$-SKT manifold, for $t\in\Delta$ sufficiently close to 0. 

On the other hand, by \cite{Wu}, the $\partial\bar{\partial}$-lemma property is open under holomorphic deformations.
\end{Conc}

So, we have proved the following:

\begin{The}\label{SKT}
Let $\pi: \mathcal{X}\longrightarrow \Delta$ be a holomorphic family of compact complex n-dimensional manifolds. Fix $p\in \{1,\cdots,n-1\}$. If  $X_0$ is a $p-$SKT $\partial\bar{\partial}$-manifold, then  $X_t$ is a $p-$SKT $\partial\bar{\partial}$-manifold for $t\in\Delta$ sufficiently close to $0$.
\end{The}

\section{Deformation limits of positive cones}
Let us start by defining the notions of $E_k$-closedness and $\overline{E_k}$-closedness.
\begin{Def}
A smooth $(r,s)$-form $\alpha^{r,s}$ is called {\bf $E_k$-closed} (resp. {\bf $\overline{E_k}$-closed}) if $\bar{\partial}\alpha^{r,s}=0$, $\partial\alpha^{r,s}=\bar{\partial}\alpha^{r+1,s-1}$, $\cdots$, $\partial\alpha^{r+k-1,s-k+1}=\bar{\partial}\alpha^{r+k,s-k}$ (resp. $\partial\alpha^{r,s}=0$, $\bar{\partial}\alpha^{r,s}=\partial\alpha^{r+1,s-1}$, $\cdots$, $\bar{\partial}\alpha^{r+k-1,s-k+1}=\partial\alpha^{r+k,s-k}$) with $\alpha^{r+l,s-l}\in C^\infty_{r+l,s-l}(X,\C)$ for $l\in\{1,\cdots,k\}$.
\end{Def}
Recall that a $C^\infty$  strictly weakly positive $(p,p)-$form $\Omega$  is $p-$Hermitian-symplectic ($p-$HS) if and only if there exist
 $\alpha^{i,2p-i}\in C_{i,2p-i}^\infty(X,\mathbb{C})$ for $i=0,\cdots,2p$ such that $$d(\sum_{i=0}^{p-1}\alpha^{i,2p-i}+\Omega+\sum_{i=0}^{p-1}\overline{\alpha^{i,2p-i}})=0.$$
 Therefore,\\

\noindent $\Omega$ is $p$-HS \, $\Longleftrightarrow$ \, 
$\partial\overline{\alpha^{0,2p}}=0,\, \partial\alpha^{2p-1,1}+\bar{\partial}\alpha^{2p,0}=0,\, \cdots, \, \partial\alpha^{p+1,p-1}+\bar{\partial}\alpha^{p+2,p-2}=0 \, \mbox{ and } \, \partial\Omega+\bar{\partial}\alpha^{p+1,p-1}=0$\\

\hspace{1.1cm} $\Longleftrightarrow$ $\alpha^{2p,0}$ is $\overline{E}_{p+1}$-closed \, and \, $\overline{d}_{p+1}(\{\alpha^{2p,0}\}_{\overline{E}_{p+1}})=\{\bar{\partial}\Omega\}_{\overline{E}_{p+1}}$\\

\hspace{1.1cm} $\Longleftrightarrow$ $\alpha^{0,2p}$ is $E_{p+1}$-closed \, and \, $d_{p+1}(\{\alpha^{0,2p}\}_{E_{p+1}})=\{\partial\Omega\}_{E_{p+1}}$.\\

\noindent So, given a strictly weakly positive $(p,p)$-form $\Omega$, we have:\\

$\Omega$ is $p$-HS \, $\Longleftrightarrow$ \, $\partial\Omega$ is $E_{p+1}$-closed and $\{\partial\Omega\}_{E_{p+1}}\in \mbox{ Im } d_{p+1}$ \, $\Longleftrightarrow$ \, $\partial\Omega$ is  $E_{p+2}$-exact.\\

\noindent This is equivalent to the property
\begin{equation}\label{E_k}
 (E_k) : \quad \bar{\partial}\alpha^{p+k,p-k}=-\partial\alpha^{p+k-1,p-k+1} \hspace*{1cm} \forall k=1,\cdots, p 
\end{equation}

We now introduce our main objects of study in this section.

\begin{Def}\label{Def:cones} Let $X$ be a compact $n$-dimensional complex manifold. Let $p\in\{1,\dots , n\}$. The cones $\mathcal{A}_p(X)$ and $\mathcal{C}_p(X)$ are defined as:
\begin{equation}\label{A_p}
\mathcal{A}_p(X)=\{[\Omega]_A/ \Omega \mbox{ strictly weakly positive such that } \partial\bar{\partial}\Omega=0\}\subset H_A^{p,p}(X,\mathbb{R})\subset H_A^{p,p}(X,\mathbb{C}),
\end{equation}
\begin{equation}\label{C_p}
\mathcal{C}_p(X)=\{[\Omega]_A/ \Omega \mbox{ strictly weakly positive such that } \Omega \mbox{ is } p\mbox{-HS} \}\subset H_A^{p,p}(X,\mathbb{R})\subset H_A^{p,p}(X,\mathbb{C}).
\end{equation}

\end{Def}

Note that $\mathcal{C}_p(X)\subset\mathcal{A}_p(X)$.

\begin{Lem}\label{Lem:open-convex-cones} The subsets $\mathcal{C}_p(X)\subset\mathcal{A}_p(X)$ are open convex cones in $ H_A^{p,p}(X,\mathbb{R})$.

\end{Lem}

\noindent {\it Proof.} To see that $\mathcal{A}_p(X)$ is a convex cone, let $[\Omega]_A,[\tilde{\Omega}]_A \in\mathcal{A}_p(X)$ and $\lambda>0$. Then $\Omega+\lambda\tilde{\Omega}$ remains strictly weakly positive and $\partial\bar{\partial}(\Omega+\lambda\tilde{\Omega})=0$, i.e. $[\Omega+\lambda\tilde{\Omega}]_A\in\mathcal{A}_p(X)$. 

To prove openness for $\mathcal{A}_p(X)$, let $[\Omega]_A \in\mathcal{A}_p(X)$ and $[\gamma]_A\in H^{p,p}_A(X,\mathbb{C})$ be arbitrary. We will prove that $[\Omega]_A+\varepsilon[\gamma]_A\in\mathcal{A}_p(X)$ for all sufficiently small $\varepsilon>0$. 

Pick arbitrary representatives $\Omega$ and $\gamma$ of their respective Aeppli cohomology classes such that $\Omega$ is strictly weakly positive. Notice that $\partial\bar{\partial}(\Omega+\varepsilon\gamma)=0$ for every $\varepsilon$ since $\Omega$ and $\gamma$ represent Aeppli classes.

It remains to prove that $\Omega+\varepsilon\gamma$ is strictly weakly positive for all sufficiently small $\varepsilon>0$. Since this positivity is a pointwise property, we can reason locally. Let $x_0\in X$ be an arbitrary point.    

By hypothesis, $(\Omega\wedge i \alpha_1\wedge \overline{\alpha_1}\wedge\cdots\wedge i \alpha_{n-p}\wedge \overline{\alpha_{n-p}})(x_0)>0$ for every choice of locally defined, linearly independent, smooth $(1,\,0)$-forms $\alpha_1,\cdots \alpha_{n-p}$. Equivalently, the restriction of $\Omega(x_0)$ to the $p$-dimensional vector subspace $E_0=E_{x_0}=\ker\alpha_1(x_0)\cap\dots\cap\ker\alpha_{n-p}(x_0)\subset T^{1,\,0}_{x_0}X$ is a {\it positive} (i.e. $>0$) volume form on $E_0=E_{x_0}$. Now, the quantity $(\Omega+\varepsilon\gamma)(x)_{|E_x}$ depends in a $C^\infty$ way on the triple $(\varepsilon,\,x,\,E)$. Therefore, by continuity, there exist a constant $\varepsilon_0>0$ and open neighbourhoods $U_{x_0}$ of $x_0$ in $X$ and ${\cal U}_{E_0}$ of $E_0$ in the Grassmannian $G_{p,\,n}$ of $p$-dimensional vector subspaces of $\C^n$, such that

$$(\Omega+\varepsilon\gamma)_{|E_x}>0 \hspace{3ex} \mbox{for all} \hspace{1ex} (\varepsilon,\, x,\,E)\in(0,\,\varepsilon_0)\times U_{x_0}\times{\cal U}_{E_0}.$$

\noindent Since both $X$ and the Grasssmannian $G_{p,\,n}$ are {\it compact} manifolds, we conclude that there exists a uniform constant $\varepsilon_0>0$ such that the $(p,\,p)$-form $\Omega+\varepsilon\gamma$ is strictly weakly positive for all $0<\varepsilon<\varepsilon_0$.

Therefore, the cone $\mathcal{A}_p(X)$ is open. The same arguments apply to $\mathcal{C}_p(X)$. \hfill $\Box$

\vspace{3ex}

We now show that the cones $\mathcal{A}_p(X)$ and $\mathcal{C}_p(X)$ are equal if the hypotheses $(H_k)$ spelt out in (\ref{H_k}) below for $k=1,\cdots, p+1$ are satisfied. These hypotheses are a collection of special cases of the $\partial\bar\partial$-lemma in a few select bidegrees. Thus, it is not necessary to assume the validity of the $\partial\bar\partial$-lemma in full generality.

\begin{Prop}
Let $X$ be a compact complex manifold with $\dim_{\mathbb{C}}X=n$. For a fixed $p\in\{1,\dots , n-1\}$ and a fixed $k\in\{1,\dots,p+1\}$, let us consider the following hypothesis $(H_k)$:

\begin{equation}\label{H_k}
 (H_k): \quad \forall\,\Gamma\in C_{p+k,p-k+1}^\infty(X,\mathbb{C}) \mbox{ such that } d\Gamma=0,\, \Gamma\in Im\,\partial\Rightarrow\Gamma\in Im\,\bar{\partial}.
  \end{equation}

\begin{enumerate}
\item[(i)] If the hypotheses $(H_1),\cdots, \mbox{ and } (H_{p+1})$ are satisfied, then $\mathcal{A}_p(X)=\mathcal{C}_p(X)$. 

\item[(ii)]If $\mathcal{A}_p(X)=\mathcal{C}_p(X)$, the hypothesis $(H_1)$ holds. 
\end{enumerate}
\end{Prop}

\noindent {\it Proof}.
\begin{enumerate}
\item[(i)]Assume that the hypotheses $(H_1),\cdots, \mbox{ and } (H_{p+1})$ are satisfied. Let $[\Omega]_A\in\mathcal{A}_p(X)$. Since $\partial\bar{\partial}\Omega=0$, $\partial\Omega$ is a $d$-closed and $\partial$-exact $(p+1,p)-$form  and by the hypothesis $(H_1)$ there exists $\alpha^{p+1,p-1}\in C_{p+1,p-1}^\infty(X,\mathbb{C})$ that solves equation $(E_1):$ $\bar{\partial}\alpha^{p+1,p-1}=-\partial\Omega$. Let $\alpha^{p+1,p-1}\in C_{p+1,p-1}^\infty(X,\mathbb{C})$ be an arbitrary solution of the equation $(E_1)$. On the other hand we have $\partial\alpha^{p+1,p-1}\in C_{p+2,p-1}^\infty(X,\mathbb{C})$. This implies that $\partial\alpha^{p+1,p-1}$ is a $d$-closed and $\partial$-exact $(p+2,p-1)$-form. Hence, under the $(H_2)$ assumption we have $\partial\alpha^{p+1,p-1}\in Im\, \bar{\partial}$, so there exists $\alpha^{p+2,p-2}\in C^\infty_{p+1,p-1}(X,\mathbb{C})$ that solves $(E_2)$  i.e. $\partial\alpha^{p+1,p-1}=-\bar{\partial}\alpha^{p+2,p-2}$.  \\
We follow the same process until we obtain, under the hypothesis $(H_p)$, the equation $(E_p)$ had a solution $\alpha^{2p,0}$ such that $\bar{\partial}\alpha^{2p,0}=-\partial\alpha^{2p-1,1}$.  This implies that $\bar{\partial}\partial\alpha^{2p,0}=0$, so $\partial\alpha^{2p,0}$ is a $(2p+1,0)$-form that is $d$-closed and $\partial$-exact. This implies that $\partial\alpha^{2p,0}\in Im\, \bar{\partial} $, hence $\partial\alpha^{2p,0}=0$.\\
 Therefore $d(\sum_{i=0}^{p-1}\alpha^{i,2p-i}+\Omega+\sum_{i=0}^{p-1}\overline{\alpha^{i,2p-i}})=0$, i.e $[\Omega]_A\in\mathcal{C}_p(X)$.\\
 It follows that $\mathcal{A}_p(X)=\mathcal{C}_p(X)$.
\item[(ii)] Consider the following linear map :
$$
\begin{array}{ll}
&H_A^{p,p}(X,\mathbb{C})\overset{T}{\longrightarrow} H_{\bar{\partial}}^{p+1,p}(X,\mathbb{C})\\
& \hspace{1.1cm}[\Omega]_A  \longmapsto [\partial\Omega]_{\bar{\partial}}
\end{array}
$$
Let us show that this map is well defined. Let $[\Omega]_A\in H_A^{p,p}(X,\mathbb{C})$,  hence $\partial\bar{\partial}\Omega=0$, i.e $\bar{\partial}(\partial\Omega)=0$, thus $\partial\Omega$ defines a class  $[\partial\Omega]_{\bar{\partial}}\in H_{\bar{\partial}}^{p,p}(X,\mathbb{C})$. If $[\Omega_1]_A=[\Omega_2]_A$, then there exists a $(p-1,p)$-form $u$ and a $(p,p-1)$-form $v$ such that $\Omega_1-\Omega_2=\partial u +\bar{\partial} v$. Thus $\partial(\Omega_1-\Omega_2)=\partial\bar{\partial}v=\bar{\partial}(-\partial v)\in Im\, \bar{\partial}$. This proves that $[\partial\Omega_1]_{\bar{\partial}}=[\partial\Omega_2]_{\bar{\partial}}$.\\
Consequently the map $T$ is well defined.

\noindent Now, we want to prove that the map $T$ vanishing identically is equivalent to the hypothesis $(H_1)$. Indeed, let $\Gamma\in C_{p+1,p}^\infty(X,\mathbb{C})$ such that $\Gamma=\partial\Omega\in \ker d$, so $\bar{\partial}\Gamma=\bar{\partial}\partial\Omega=-\partial\bar{\partial}\Omega$. But $\Gamma$ is a $d$-closed form of pure type $(p+1,p)$, i.e $\partial\Gamma=0$ and $\bar{\partial}\Gamma=0$, then $\partial\bar{\partial}\Omega=0$.

On the other hand $T\equiv 0$ means that $T([\Omega]_A)=[\partial\Omega]_{\bar{\partial}}=0\in H_{\bar{\partial}}^{p+1,p}(X,\mathbb{C})$, for every $[\Omega]_A\in H^{p,p}_A(X,\C)$. Hence $\partial\Omega\in Im\, \bar{\partial}$.

Conversely, assume that for all  $\Gamma\in C_{p+1,p}^\infty(X,\mathbb{C})$ such that $d\Gamma=0$, we have the implication  
$\Gamma\in Im \,\partial$ $\Rightarrow$ $\Gamma\in Im\,\bar{\partial}$. Let $[\Omega]_A\in H^{p,p}_A(X,\C)$. Then $\Gamma=\partial\Omega$  implies that $\Gamma=\partial\Omega\in Im\,\bar{\partial}$, so $[\partial\Omega]_{\bar{\partial}}=0=T([\Omega]_A)$. This proves  $T\equiv 0$.
 
Now, since $\mathcal{A}_p(X)=\mathcal{C}_p(X)$, we have  $$\mathcal{A}_p(X)\,\cap\,\ker T=\{[\Omega]_A/ \Omega \mbox{ strictly weakly positive such that } \partial\bar{\partial}\Omega=0 \mbox{ and } \partial\Omega\in Im\, \bar{\partial}\}\supset \mathcal{C}_p(X).$$ 
So $\mathcal{C}_p(X)\subset \mathcal{A}_p(X)\,\cap\,\ker T\subset \mathcal{A}_p(X)$. Since $\mathcal{A}_p(X)=\mathcal{C}_p(X)$, this implies that  $\mathcal{A}_p(X)\cap\ker T=\mathcal{A}_p(X)$.\\
Hence  $\ker T=H_A^{p,p}(X,\mathbb{C})$, which is equivalent to the $(H_1)$ assumption.
\end{enumerate} \hfill $\Box$

\noindent We shall now consider, for $k=1,\cdots,p+1$, the linear map:
$$
\begin{array}{ll}
&H_A^{p+k-1,p-k+1}(X,\mathbb{C})\overset{\widehat T_k}{\longrightarrow} H_{BC}^{p+k,p-k+1}(X,\mathbb{C})\\
&\hspace*{2.7cm}[\Omega]_A  \longrightarrow [\partial\Omega]_{BC}
\end{array}
$$
This map is well defined. Indeed.\\
Suppose that $\partial\bar{\partial}\Omega=0$, i.e. $\bar{\partial}(\partial\Omega)=0$, so $\partial\Omega\in\ker \bar{\partial}$ this is equivalent to $\partial\Omega\in \ker\partial\cap\ker\bar{\partial}$.  Moreover, if $\Omega=\partial u +\bar{\partial} v$  with  $u$ and $v$  are $(p+k-2,p-k+1)$-form and $(p+k-1,p-k)$-form respectively, hence $\partial\Omega=\partial\bar{\partial} v\in Im\, \partial\bar{\partial}$. Thus the map $\widehat T_k$ is well defined.\\

Consider the linear map :
$$
\begin{array}{ll}
I^{p+k,p-k+1} :&H_{BC}^{p+k,p-k+1}(X,\mathbb{C})\longrightarrow H_{\partial}^{p+k,p-k+1}(X,\mathbb{C})\\
&\hspace*{2cm}[\Gamma]_{BC}  \longrightarrow [\Gamma]_{\partial}
\end{array}
$$
We have $\partial\Gamma=0$, and $\bar{\partial}\Gamma=0$. To show that the map $I^{p+k,p-k+1}$ is well defined, we still have to show that the definition is independent of the choice of representative of the class $[\Gamma]_{BC} $. In other words, if $\Gamma=\partial\bar{\partial}u$  with $u$ a $(p+k-1,p-k)$-form, then $\Gamma\in Im\, \partial$. This is obvious, so the map $I^{p+k,p-k+1}$ is well defined.\\

\noindent We will infer the following
\begin{Rem}
The map $\widehat{T}_k$ vanishes identically if and only if for all $ \Omega\in C_{p+k-1,p-k+1}^\infty (X,\mathbb{C})\cap\ker\partial\bar{\partial}$, we have $ \partial\Omega\in Im\, \partial\bar{\partial}$. This is equivalent to:
\begin{equation}\label{tilde H_k}
(\tilde{H}_k): \quad \forall\Gamma\in C_{p+k,p-k+1}^\infty (X,\mathbb{C}) \mbox{ such that } d\Gamma=0, \mbox{ we have } \Gamma\in Im\, \partial \Rightarrow \Gamma\in Im\, \partial\bar{\partial}.
\end{equation}

 \noindent This is further equivalent to $I^{p+k,p-k+1}$ is injective, hence to $\ker I^{p+k,p-k+1}=\{0\}$.\\
 
 Consider the diagram:
$$
\xymatrix{
H_A^{p+k-1,p-k+1}(X,\mathbb{C}) \ar[dr]_{g} \ar[r]^{\widehat{T}_k} & H_{BC}^{p+k,p-k+1}(X,\mathbb{C}) \ar[d]^{I^{p+k,p-k+1}}\\
                                               & H_{\partial}^{p+k,p-k+1}(X,\mathbb{C}) 
                                               }
$$
where $g=I^{p+k,p-k+1}\circ  \widehat{T}_k $. We have $\ker I^{p+k,p-k+1}=\{[\Gamma]_{BC}/[\Gamma]_\partial=0\}  $  and $Im \, \widehat{T}_k=\{[\partial\Omega]_{BC}/ \Omega\in\ker\partial\bar{\partial}\subset C_{p+k-1,p-k+1}^\infty (X,\mathbb{C})\} \subset \ker I^{p+k,p-k+1} $.

Conversely,  if $[\Gamma]_{BC}\in \ker I^{p+k,p-k+1}$, then there exists a $(p+k-1,p-k+1)$-form $\Omega$ such that $\Gamma=\partial\Omega$. Since $[\Gamma]_{BC}\in H_{BC}^{p+k,p-k+1}(X,\mathbb{C})$, then $0=\bar{\partial}\Gamma=\bar{\partial}\partial\Omega$. This implies that $[\Gamma]_{BC}=[\partial\Omega]_{BC}=\widehat{T}_k ([\Omega]_A)$. Therefore we always have:  $$\ker I^{p+k,p-k+1}=Im \,\widehat{T}_k$$
Furthermore, the hypothesis $(\tilde{H}_k)$ (cf. \ref{tilde H_k}) is satisfied if and only if $\ker I^{p+k,p-k+1}=Im\, \widehat{T}_k=\{0\}$.
\end{Rem}

\noindent We can now prove the following
\begin{The}
Let $(X_t)_{t\in\Delta}$ be a holomorphic family of compact complex manifolds with  $\dim_{\mathbb{C}} X_t=n$ and $(\omega_t)_{t\in\Delta}$ a smooth family of Hermitian metrics on $(X_t)_{t\in\Delta}$ for $t\in\Delta$. If for all $t$ sufficiently close to $0$ we have :
$$
 \left\{
\begin{array}{rl}
& h_A^{p+k-1,p-k+1}(0)=h_A^{p+k-1,p-k+1}(t):=dim_\mathbb{C} H_A^{p+k-1,p-k+1}(X_t,\mathbb{C}) \\
&\\
& h_{BC}^{p+k,p-k+1}(0)=h_{BC}^{p+k,p-k+1}(t):=dim_\mathbb{C} H_{BC}^{p+k,p-k+1}(X_t,\mathbb{C}) \\
&\\
& h_\partial^{p+k,p-k+1}(0)=h_\partial^{p+k,p-k+1}(t):=dim_\mathbb{C} H_{\partial}^{p+k,p-k+1}(X_t,\mathbb{C}),
\end{array}
\right.
$$
$$
\mbox{ then \hspace*{0.5cm} }\left\{
\begin{array}{rl}
&\Delta\ni t \overset{\mathcal{H}_A}{\longmapsto} H_A^{p+k-1,p-k+1}(X_t,\mathbb{C})  \\
&\\
&\Delta\ni t \overset{\mathcal{H}_{BC}}{\longmapsto} H_{BC}^{p+k,p-k+1}(X_t,\mathbb{C})\\
&\\
&\Delta\ni t \overset{\mathcal{H}_{\partial}}{\longmapsto} H_{\partial}^{p+k,p-k+1}(X_t,\mathbb{C})
\end{array}
\right.
$$
are $C^\infty$ vector bundles. Moreover, the linear maps:
$$
\xymatrix{
H_A^{p+k-1,p-k+1}(X_t,\mathbb{C}) \ar[dr]_{g_k(t)} \ar[r]^{\widehat{T}_k(t)} & H_{BC}^{p+k,p-k+1}(X_t,\mathbb{C}) \ar[d]^{I^{p+k,p-k+1}(t)}\\
                                               & H_{\partial_t}^{p+k,p-k+1}(X_t,\mathbb{C}) 
                                               }
$$
vary in a $C^\infty$ way with $t\in\Delta$, where $g_k=I^{p+k,p-k+1}\circ  \widehat{T}_k $.
\end{The}

\noindent \noindent {\it Proof}. The Laplace-type operators
$$
\left\{
\begin{array}{rl}
&\Delta_{A,t}=\partial_t\partial_t^*+\bar{\partial}_t\bar{\partial}_t^*+\bar{\partial}_t^*\partial_t^*\partial_t\bar{\partial}_t+\partial_t\bar{\partial}_t\bar{\partial}_t^*\partial_t^*+\partial_t\bar{\partial}_t^*\bar{\partial}_t\partial_t^*+\bar{\partial}_t\partial_t^*\partial_t\bar{\partial}_t^* \hspace{2cm} (cf.\quad \cite{Sch07})\\
&\Delta_{BC,t}=\partial_t\bar{\partial}_t\bar{\partial}_t^*\partial_t^*+\bar{\partial}_t^*\partial_t^*\partial_t\bar{\partial}_t+\bar{\partial}_t^*\partial_t\partial_t^*\bar{\partial}_t+\partial_t^*\bar{\partial}_t\bar{\partial}_t^*\partial_t+\bar{\partial}_t^*\bar{\partial}_t+\partial_t^*\partial_t \hspace{2cm} (cf.\quad \cite{KS60})\\
&\Delta_t'=\partial_t\partial_t^*+\partial_t^*\partial_t
\end{array}
\right.
$$
are elliptic (cf. \cite{KS60}, \cite{Sch07}), so we have the Hodge isomorphisms:
$$
\left\{
\begin{array}{rl}
&H_A^{p+k-1,p-k+1}(X_t,\mathbb{C}) \simeq \mathcal{H}_{\Delta_{A,t}}^{p+k-1,p-k+1}(X_t,\mathbb{C})\\
&\\
&H_{BC}^{p+k,p-k+1}(X_t,\mathbb{C}) \simeq \mathcal{H}_{\Delta_{BC,t}}^{p+k,p-k+1}(X_t,\mathbb{C})\\
&\\
&H_{\partial_t}^{p+k,p-k+1}(X_t,\mathbb{C}) \simeq \mathcal{H}_{\Delta'_t}^{p+k,p-k+1}(X_t,\mathbb{C}).
\end{array}
\right.
$$
where  $\mathcal{H}_{\Delta_{BC,t}}^{p+k,p-k+1}(X_t,\mathbb{C})=\ker\Delta_{BC,t}$ and $\mathcal{H}_{\Delta_{A,t}}^{p+k-1,p-k+1}(X_t,\mathbb{C})=\ker\Delta_{A,t}$ stand for the Bott-Chern and Aeppli harmonic spaces (spaces of Bott-Chern and Aeppli harmonic forms ) respectively and $ \mathcal{H}_{\Delta'_t}^{p+k,p-k+1}(X_t,\mathbb{C})=\ker \Delta'_t$.\\
On the otherhand, since $h_A^{p+k-1,p-k+1}(0)=h_A^{p+k-1,p-k+1}(t)$, $ h_{BC}^{p+k,p-k+1}(0)=h_{BC}^{p+k,p-k+1}(t)$ and  $h_\partial^{p+k,p-k+1}(0)=h_\partial^{p+k,p-k+1}(t)$ for all $t$ sufficiently close to $0$, by Kodaira-spencer   (cf. \cite{K86}, theorem 7.4) we have :
 $$
\left\{
\begin{array}{rl}
& t \overset{\mathcal{H}_A}{\longmapsto} H_A^{p+k-1,p-k+1}(X_t,\mathbb{C}) \simeq \mathcal{H}_{\Delta_{A,t}}^{p+k-1,p-k+1}(X_t,\mathbb{C}) \\
&\\
& t \overset{\mathcal{H}_{BC}}{\longmapsto} H_{BC}^{p+k,p-k+1}(X_t,\mathbb{C})\simeq \mathcal{H}_{\Delta_{BC,t}}^{p+k,p-k+1}(X_t,\mathbb{C})\\
&\\
& t \overset{\mathcal{H}_\partial}{\longmapsto} H_{\partial}^{p+k,p-k+1}(X_t,\mathbb{C})\simeq \mathcal{H}_{\Delta'_t}^{p+k,p-k+1}(X_t,\mathbb{C})
\end{array}
\right.
$$
are $C^\infty$ vector bundles. Recall that the linear maps $\widehat{T}_k(t)$ and $I^{p+k,p-k+1}(t)$ are well defined as shown above. Let $h_t$ and $F_t$ be the orthogonal projection of $C_{p+k-1,p-k+1}^\infty(X_t)$, resp. $C_{p+k,p-k+1}^\infty(X_t)$ onto $\mathcal{H}_A^{p+k-1,p-k+1}(X_t,\mathbb{C})$, resp. $\mathcal{H}_{BC}^{p+k,p-k+1}(X_t,\mathbb{C}) $ respectively.
$$
\xymatrix{
C_{p+k-1,p-k+1}^\infty(X_t) \ar[d]_{h_t} \ar[r]^{f_t=\partial_t} & C_{p+k,p-k+1}^\infty(X_t) \ar[d]^{F_t} \\
H_A^{p+k-1,p-k+1}(X_t,\mathbb{C})\simeq\mathcal{H}_{\Delta_A}^{p+k-1,p-k+1}(X_t,\mathbb{C}) \ar[r]_{\widehat{T}_k(t)}       & \mathcal{H}_{\Delta_{BC}}^{p+k,p-k+1}(X_t,\mathbb{C}) \simeq H_{BC}^{p+k,p-k+1}(X_t,\mathbb{C})
}
$$
We have $F_t\circ f_t=\widehat{T}_k(t)\circ h_t$ , $f_t$ is $C^\infty$ and  by Kodaira-Spencer as in \cite{K86} $F_t$ and $h_t$ vary smoothly with $t\in\Delta$. Therefore $\widehat{T}_k(t)$ varies smoothly with $t$. Let $P_t$ the orthogonal projection of $C_{p+k,p-k+1}^\infty(X_t)$ onto $\mathcal{H}_{\Delta'_t}^{p+k,p-k+1}(X_t,\mathbb{C}) $.
$$
\xymatrix{
C_{p+k,p-k+1}^\infty(X_t) \ar[dr]^{P_t} \ar[d]_{F_t}        \\
\mathcal{H}_{\Delta_{BC}}^{p+k,p-k+1}(X_t,\mathbb{C})\quad  \ar[r]_{I^{p+k,p-k+1}(t)}  &   \quad \mathcal{H}_{\Delta'_t}^{p+k,p-k+1}(X_t,\mathbb{C}) 
}
$$
We have $P_t =I^{p+k,p-k+1}(t)\circ F_t$ and by Kodaira-Spencer as in \cite{K86} $P_t$ and $F_t$ are $C^\infty$ with $t\in\Delta$, then $I^{p+k,p-k+1}(t)$ is also $C^\infty$ with $t\in\Delta$.
Consequently $g(t)$ varies smoothly with $t\in\Delta$. \hfill $\Box$

\vspace{2ex}

One obtains sections $\widehat{T}_k=(\widehat{T}_k(t))_{t\in\Delta}\in C^\infty(\Delta,End(\mathcal{H}_A,\mathcal{H}_{BC}))$ and $I^{p+k,p-k+1}=(I^{p+k,p-k+1}(t))_{t\in\Delta}\in C^\infty(\Delta,End(\mathcal{H}_{BC},\mathcal{H}_{\partial}))$.

\noindent As a consequence, we obtain the following conclusion on the deformation limit of the $X_t$'s that satisfy the $(\tilde{H}_k)$ assumption.

\begin{Cor}
Let $(X_t)_{t\in\Delta}$ be a holomorphic family of compact complex manifolds and $(\omega_t)_{t\in\Delta}$ a smooth family of metrics on $(X_t)_{t\in\Delta}$, $\dim_{\mathbb{C}} X_t=n$, $t\in\Delta$.
$$
\mbox{Assume that :\hspace*{0.5cm} }\left\{
\begin{array}{rl}
& h_A^{p+k-1,p-k+1}(0)=h_A^{p+k-1,p-k+1}(t)\\
&\\
& h_{BC}^{p+k,p-k+1}(0)=h_{BC}^{p+k,p-k+1}(t) \hspace*{2cm} \forall t\sim 0\\
&\\
& h_\partial^{p+k,p-k+1}(0)=h_\partial^{p+k,p-k+1}(t)
\end{array}
\right.
$$
If $X_t$ satisfies the hypothesis $(\tilde{H}_k)$, $\forall t\in \Delta\setminus\{0\}$, then $X_0$ also satisfies the hypothesis $(\tilde{H}_k)$ .
\end{Cor}

\noindent \noindent {\it Proof}. Recall that $X_t$ satisfying $(\tilde{H}_k)$ is equivalent to the map $\widehat{T}_k(t)$ vanishing identically.

\noindent Moreover, if $\widehat{T}_k(t)\equiv 0$, for all $t\in \Delta\setminus\{0\}$, then by continuity of $\widehat{T}_k(t)$, $\widehat{T}_k(0)\equiv 0$ which is equivalent to $X_0$ satisfying the $(\tilde{H}_k)$ assumption. \hfill $\Box$

\vspace{2ex}

Hence we obtain the following:
\begin{Prop}\label{Prop} Fix $p\in\{1,\cdots,n-1\}$ and let $k\in\{1,\cdots,p+1\}$ such that $p+k\leq n$.
\begin{enumerate}
\item[(1)] $\forall k\in\{1,\cdots,p+1\}$, $(\tilde{H}_k)$ $\Rightarrow$ $(H_k)$\\
\item[(2)] $(\tilde{H}_1)+\cdots+(\tilde{H}_{p+1})$ $\Rightarrow$ $(H_1)+\cdots+(H_{p+1})$ $\Rightarrow$ $\mathcal{A}_p(X)=\mathcal{C}_p(X)$\\
\item[(3)] Suppose that $\mathcal{A}_p(X_t)=\mathcal{C}_p(X_t)$  (cf. (\ref{A_p}),(\ref{C_p})) $\forall t\in\Delta\setminus\{0\} $. Then:
$$
\left.
\begin{array}{rl}
&\forall k\in\{1,\cdots,p+1\}\\
& (\tilde{H}_1)+\cdots+(\tilde{H}_{p+1}) \hspace*{2cm}\forall t\in\Delta\setminus\{0\}\\
&\\
& h_A^{p+k-1,p-k+1}(0)=h_A^{p+k-1,p-k+1}(t)\\
&\\
& h_{BC}^{p+k,p-k+1}(0)=h_{BC}^{p+k,p-k+1}(t) \hspace*{1cm} \forall t\sim 0\\
&\\
& h_\partial^{p+k,p-k+1}(0)=h_\partial^{p+k,p-k+1}(t)
\end{array}
\right\}\Rightarrow \mathcal{A}_p(X_0)=\mathcal{C}_p(X_0)
$$
\end{enumerate}
\end{Prop}

\begin{Rem}
In the case where $n=3$ and $p=2$, $k\in\{1,2,3\}$ and $2+k\leq 3 $ we must have $k=1$.\\
So, $\mathcal{A}_2(X)=\{[\Omega]_A/\Omega>0 \mbox{ such that } \partial\bar{\partial}\Omega=0\}=\mathcal{G}_X$
 Gauduchon cone\\
$\mathcal{C}_2(X)=\{[\Omega]_A/\Omega>0 \mbox{ and there exists }\alpha^{1,3}\in C_{1,3}^\infty (X,\mathbb{C}) \mbox{ such that } d(\alpha^{1,3}+\Omega+\overline{\alpha^{1,3}})=0\}=\mathcal{SG}_X$  sG cone and the following equivalence holds:\\
 $$\mathcal{A}_2(X)=\mathcal{C}_2(X) \Leftrightarrow X \mbox{ is sGG i.e } \mathcal{SG}_X=\mathcal{G}_X$$
\end{Rem}

\begin{Prop}
Let $p,q\in\{1,\cdots,n\}$ be fixed. Suppose that the implication :
$$ u\in \mbox{Im }\, \partial  \quad \Rightarrow \quad u\in \mbox{Im }\, \partial\bar{\partial} $$ holds for all $d$-closed forms of types $(p,q)$, $(q,p)$, $(p+1,q)$, and $(q+1,p)$ for all $p,q$ such that $p+q=k$. Then, there exists a canonical injective linear map :
$$
\begin{array}{ll}
& H_A^{p,q}(X,\mathbb{C})\hookrightarrow H_{DR}^{p+q}(X,\mathbb{C})\\
&\hspace{1.1cm}[\alpha]_A\mapsto \{\alpha\}_{DR}
\end{array}
$$
where $\alpha$ is any $d$-closed representative of the class $[\alpha]_A$ (such an $\alpha$ exists due to the hypothesis)
\end{Prop}

\noindent \noindent {\it Proof}. The map
$$
\begin{array}{ll}
& H_A^{p,q}(X,\mathbb{C})\longrightarrow H_{DR}^{p+q}(X,\mathbb{C})\\
&\hspace{1.1cm}[\alpha]_A\mapsto \{\alpha\}_{DR}
\end{array}
$$
is well defined since if $\partial\bar{\partial}\alpha=0$, then by the hypothesis $Im\, \partial \subset Im\, \partial\bar{\partial}$, there exists $(p,q-1)$-form $v$ such that $\partial\alpha=-\partial\bar{\partial}v$. On the other hand, $\bar{\alpha}$ is a $(q,p)$-form and $\partial\bar{\alpha}$ is a $\partial$-exact $(q+1,p)$-form. Hence, by assumption, $\partial\bar{\alpha}$ is $\partial\bar{\partial}$-exact. Then, by conjugation, $\bar{\partial}\alpha$ stills $\partial\bar{\partial}$-exact, i.e. there exists some $(p-1,q)$-form $u$ such that $\bar{\partial}\alpha=\partial\bar{\partial}u$. This implies that $d(\alpha+\partial u +\bar{\partial} v)=0$, it means that every Aeppli cohomology class contains a $d$-closed representative.

Let $\alpha$ to be a $(p,q)$-form such that $d\alpha=0$ and $\alpha=\partial u + \bar{\partial} v$. Then  $\partial\alpha=0=\partial\bar{\partial}v$ and $\bar{\partial}\alpha=0=\bar{\partial}\partial u$. Note that $\partial u$ is a $d$-closed $\partial$-exact $(p,q)$-form, so $\partial u\in Im\, \partial\bar{\partial}$. Thus $\partial u\in Im\, d$. Meanwhile,
$\bar{\partial} v$ is a $d$-closed $\bar{\partial}$-exact $(p,q)$-form, so $\bar{\partial}v\in Im\, \partial\bar{\partial}$, then $\bar{\partial}v\in Im\, d$. Therefore $\alpha\in Im\, d$, thus the map above is well defined.\\

\noindent {\it{Injectivity:}} for all $\alpha \in C_{p,q}^\infty(X,\mathbb{C})$ such that $\partial\bar{\partial}\alpha=0$ and $\alpha\in Im\, d$, we have $\alpha\in Im\, \partial+ Im \,\bar{\partial}$. Hence the map is injective once it is well defined. \hfill $\Box$

\vspace{1ex}
\noindent We can now infer the following

\begin{Cor} \label{Cor}
\begin{enumerate}
\item[(1)] Fix $k\in\{0,1,\cdots,2n\}$, suppose that the implication  $$u\in \mbox{Im } \partial \quad \Rightarrow \quad u\in \mbox{Im } \partial\bar{\partial}$$ holds for all $d$-closed forms of types $(p,q)$, $(q,p)$, $(p+1,q)$, and $(q+1,p)$ for all $p,q$ such that $p+q=k$.\\
Then there exists a canonical injection: $$\bigoplus_{p+q=k} H_A^{p,q}(X,\mathbb{C})\hookrightarrow H_{DR}^k(X,\mathbb{C})$$
\item[(2)] Fix $k\in\{0,1,\cdots,2n\}$. Suppose that $(*_k)$ holds, where
\begin{equation}\label{star}
\begin{array}{ll}
(*_k) :\quad &\mbox{ the implication } u\in Im\, \partial  \Rightarrow u\in Im\, \partial\bar{\partial} \mbox{ holds for all } d-closed \mbox{ forms of types }\\
& (p,q), (q,p), (p+1,q), \mbox{ and } (q+1,p) \\
&\mbox{ for all } p,q \mbox{ such that } p+q=k \mbox{ or } p+q=2n-k
\end{array}
\end{equation}
Then there exists a canonical injection :
$$ \bigoplus_{p+q=k}H_A^{p,q}(X,\mathbb{C})\oplus\bigoplus_{p+q=2n-k}H_A^{p,q}(X,\mathbb{C})\hookrightarrow H_{DR}^k(X,\mathbb{C})\oplus H_{DR}^{2n-k}(X,\mathbb{C}) $$
Due to Angella-Tomassini \cite{AT13}, we have :
$$2b_k\leq \sum_{p+q=k}h^{p,q}_A+\sum_{p+q=2n-k}h^{p,q}_A  $$
 this map is an isomorphism and $2b_k=\sum_{p+q=k}h_A^{p,q}+\sum_{p+q=2n-k}h_A^{p,q}$.
\end{enumerate}
\end{Cor}

\noindent Another immediate consequence is the following degeneration at $E_1$ of the Fr\"olicher spectral sequence.

\begin{Cor}
If $(*_k)$ (cf. (\ref{star})) assumption is satisfied for all $k\in\{0,1,\cdots,2n\}$, then $E_1(X)=E_\infty(X) $ (i.e. the Fr\"olicher spectral sequence degenerate at $E_1$).
\end{Cor}

\noindent In the end, the above results (cf. Corollary \ref{Cor}) leads to the following
\begin{Prop}
It is clear that $(*_{2p})$ (cf. (\ref{star})) implies the $(\tilde{H}_k)$ (cf. (\ref{tilde H_k}))  assumption for all $k\in\{1,\cdots,p+1\}$. If $X_0$ satisfies $(*_{2p})$ and $(*_{2p+1})$, then:
$$
\left\{
\begin{array}{rl}
&\forall k\in\{1,\cdots,p+1\}\\
& h_A^{p+k-1,p-k+1}(0)=h_A^{p+k-1,p-k+1}(t)\\
&\\
& h_{BC}^{p+k,p-k+1}(0)=h_{BC}^{p+k,p-k+1}(t) \hspace*{1cm} \forall t\sim 0\\
&\\
& h_\partial^{p+k,p-k+1}(0)=h_\partial^{p+k,p-k+1}(t)
\end{array}
\right.
$$
\end{Prop}

\vspace{1ex}

\vspace{6ex}

\noindent Ibn Tofail University, Faculty of Sciences, Departement of Mathematics PO 242 Kenitra, Morocco

\noindent Email: houda.bellitir@uit.ac.ma

\end{document}